\documentclass[letterpaper,12pt,twoside]{article}

	\usepackage[utf8]{inputenc}
	\usepackage{amsmath,amssymb,amsthm}
	\usepackage{bm}
       	\usepackage{pifont} 
        \usepackage{lmodern,fancyhdr,lastpage,nccfoots,caption}
	\usepackage{enumerate}
	\usepackage{graphicx,xcolor}

        \usepackage[all,2cell]{xy}

\usepackage[pdftex,hyperindex]{hyperref} 
\usepackage{placeins}

\theoremstyle{definition}

\theoremstyle{remark}

\numberwithin{equation}{section}

\usepackage[all,2cell]{xy}

		\DeclareMathAlphabet{\mathpzc}{OT1}{pzc}{m}{it}

		\renewcommand{\leq}{\leqslant}

		\theoremstyle{plain}
		\newtheorem*{non-Th}{Teorema}        
		\newtheorem*{non-Cor}{Corolario}     

		\theoremstyle{definition}

		\theoremstyle{remark}

		\numberwithin{equation}{section}

		\renewcommand{\phi}{\varphi}

\newcounter{a}
\ifodd\thea\else\stepcounter{a}\fi

\fancypagestyle{plain}{%
\fancyhf{}

}

\begin{document}

\thispagestyle{plain}

\begin{center}
\Large
\textsc{Introducci\'on a la teor\'ia de categor\'{i}as, lógica y topos elementales\\ para una mente curiosa}
\end{center}

\begin{center}
    25/Nov/2022
\end{center}

\begin{center}
  \textit{Jes\'{u}s E. S\'{a}nchez--Guevara\footnote{\tiny Con el apoyo de la Escuela de Matem\'{a}tica de la Universidad de Costa Rica, {\sc em}at--{\sc ucr}.}}\\
  \small Escuela de Matem\'atica, Universidad de Costa Rica {\sc em}at--{\sc ucr}\\
  \small San Jos\'e 11501, Costa Rica\\
  \small e-mail: \texttt{jesus.sanchez\_g@ucr.ac.cr}\\
  \small {\sc orcid} \href{https://orcid.org/0000-0001-8993-1538}{\tt \color{black!29!blue} 0000-0001-8993-1538}
\end{center}

\begin{center}
  \textit{Ronald A. Z\'u\~niga-Rojas\footnote{\tiny Con el apoyo de la Escuela de Matem\'atica de la Universidad de Costa Rica, {\sc em}at--{\sc ucr}, espec\'{i}ficamente a trav\'{e}s del {\sc cimpa}, mediante el proyecto {\tt 821-C1-010}.}}\\
  \small Centro de Investigaci\'on en Matem\'atica Pura y Aplicada {\sc cimpa}\\
  \small Escuela de Matem\'atica, Universidad de Costa Rica {\sc em}at--{\sc ucr}\\
  \small San Jos\'e 11501, Costa Rica\\
  \small e-mail: \texttt{ronald.zunigarojas@ucr.ac.cr}\\
  \small {\sc orcid} \href{https://orcid.org/0000-0003-3402-2526}{\tt \color{black!29!blue} 0000-0003-3402-2526}
\end{center}

\vspace{2ex}

\noindent {\sc Abstract.} 
This paper presents a study of how the theory
of categories leads to the creation of non classical logical systems.
In particular, the case of the elementary topos of
graphs, where there are three other truth values different from false
and true. The approach in this article to the theory of categories
avoids specialized mathematical training to understand it,
since it seeks to make accessible the main ideas of
this branch of mathematics to other disciplines of knowledge.
This work was presented by the second author, under the title of
``Logic and Categories'', at the I Colloquium on Logic, Epistemology
and Methodology organized by the School of Philosophy of the University
sity of Costa Rica.
\vspace{0.5cm}

\noindent {\sc Resumen.} 
En este trabajo se presenta un estudio de cómo la teoría de categorías
lleva a la creación de sistemas lógicos diferentes a los clásicos. 
En particular, se describe el caso del topos elemental de grafos, donde existen
otros tres estados de verdad diferentes al falso y verdadero.
El abordaje en este artículo de la teoría de categorías 
no necesita una formación matemática especializada para su comprensión,
ya que se busca hacer accesibles las principales ideas de esta rama de las matemáticas
a otras disciplinas del conocimiento.
Este trabajo fue presentado por el segundo autor, 
bajo el t\'{i}tulo de ``L\'{o}gica y Categor\'{i}as'', 
en el I Coloquio de L\'{o}gica, Epistemolog\'{i}a y Metodolog\'{i}a
organizado por la Escuela de Filosof\'{i}a de la Universidad de Costa Rica.

\vspace{2ex}

\noindent{\bf Palabras clave:} 
categor\'{i}as, topos elementales, l\'{o}gica, filosof\'{i}a.

\noindent{\bf MSC~2020:} Primaria {\tt 00A30}; Secundarias {\tt 03A10, 03G30}.

\tableofcontents

\section{Introducci\'{o}n}
\label{intro}

El comportamiento de las fuerzas físicas en el plano y el espacio es modelado en matemáticas con espacios vectoriales y las transformaciones lineales que operan entre ellos, lo cual imprime una gran importancia al estudio de las propiedades entre estos dos conceptos matemáticos. En \cite{saunders-general-theory-of-natural-equivalences} Saunders MacLane nos dice que estos nexos también se manifiestan en situaciones de diferente naturaleza, es decir, veríamos un comportamiento análogo con los grupos y sus homomorfismos, con los espacios topológicos y las aplicaciones continuas, con los complejos simpliciales y las transformaciones simpliciales, con los conjuntos ordenados y los mapeos que preservan el orden. Así, para abordar la generalidad de estas relaciones, sin tomar en cuenta la naturaleza particular de cada situación, se recurre al concepto matemático de categoría.

La teoría de categorías se encuentra presente de forma transversal en las matemáticas modernas. Su versatilidad la ha justificado como herramienta unificadora y capáz de describir y manipular en total generalidad propiedades intrínsecas de diferentes entidades matemáticas. En el tercer capítulo de \cite{james1999history} titulado \emph{Desarrollo del concepto de homotopía}, Ria Vanden Eynde nos indica que la preferencia de la topología por las teorías con puntos de base tiene como explicación los conceptos de categorías y funtores presentados por Eilenberg y MacLane en los años cuarenta, ya que, al abordar su estudio desde el punto de vista categórico, entre otras cosas, se pueden identificar en ellas objetos que son al mismo tiempo iniciales y finales, lo cual crea un paralelismo con categorías de propiedades similares como la categoría de anillos o la de grupos, es decir, son más adecuadas para la aplicación de métodos algebraicos. Lo cual no se sería posible si se trabajara en categorías de espacios sin puntos de base o en la categoría de los conjuntos.

El alcance de la nueva herramienta matemática de las categorías fue más allá de la reorganización de la información algebraica para clasificar, describir y hallar relaciones importantes. \cite{Marquis2011-MARTHO-8} Marquis y Reyes (2011), indican sobre Lawvere que sus trabajos en lógica simbólica lo llevaron a descubrir una forma de interpretar proposiciones lógicas mátemáticas en objetos y flechas de ciertas categorías especiales, llamadas topos elementales. Las cuales se pueden ver como entes que generalizan la categoría de los conjuntos y cuyas principales características se derivan de los topos de Grothendieck, introducidos en la década de los sesenta \cite{bourbaki2006theorie}.

Este nuevo punto de vista abrió la posibilidad de replantear las bases axiomáticas de las matemáticas, lo cual implicaría cambiar a los conjuntos como puntos de partida para la fundamentación de todas las matemáticas, ya que algunos especialistas como Jean Benabou \cite{benabou-youtube} y Alain Prouté \cite{alain-logique-cat}, indican indican que los axiomas conjuntistas se muestran insuficientes para poder describir la complejidad total de la maquinaria matemática moderna.  


Prouté \cite{alain-logique-cat}, dice que en el estudio de las categorías abordadas por Lawvere se ve que la estructura de un topos elemental es lo suficientemente rica para que se pueda considerar como un universo matemático completo.
En ellos se pueden encontrar modelos donde se pueden satisfacer o no el principio del tercero excluido o el axioma de elección. También, la lógica interna de un topos se puede abordar de forma intuicionista o constructivista, es decir, que la verdad se entiende solamente a través de una prueba y la existencia solo se reconoce a partir de una construcción explícita.

Este artículo está configurado en dos partes, la primera donde se discuten las principales características de las categorías en general y la segunda parte donde se abordan las propiedades de los topos elementales. Nuestro enfoque global en este texto es el de presentar un abordaje claro de la teoría de categorías para personas no especialistas en matemáticas con el fin de potenciar sus aportes a otras ramas cognitivas, como la filosofía.

\section{Categor\'ias}

	Una categoría se puede describir como cualquier sistema donde sea posible identificar flechas y los lugares entre los cuales se mueven. Además, la características de estas flechas van a estar controladas por algunas condiciones sobre su comportamiento. Como primer ejemplo trataremos con una categoría hecha a partir de un conjunto de personas, la cual llamaremos $\mathcal{P}$.

	En $\mathcal{P}$ vamos a tener una cantidad definida de personas, las cuales serán vistas como objetos entre los cuales vamos a trazar algunas flechas. Así, si $A$ y $B$ son dos personas de nuestra categoría, entonces habrá una flecha que parte de $A$ hacia $B$, si $A$ puede provocar una reacción de asombro en $B$. En este trazado de flechas podemos remarcar que no se excluye la posibilidad de que $A$ y $B$ sean la misma persona. De hecho, para que nuestro ejemplo funcione como categoría, sólo resta aceptar que cualquier persona siempre podrá provocar una reacción de asombro en ella misma.
	
	Veamos ahora las principales propiedades de las flechas de $\mathcal{P}$. Primero, siempre vamos a tener al menos tantas flechas como personas, pues por cada persona hay una flecha que sale y llega a ella. Segundo, si de $A$ a $B$ hay una flecha $f$, y de $B$ a otra persona $C$ hay una flecha $g$, entonces necesariamente hay una flecha $h$ de $A$ a $C$, pues si $A$ puede asombrar a $B$ y $B$ puede asombrar a $C$,  entonces es correcto decir que $A$ puede provocar asombro en $C$, ya que $A$ podría hacerlo al usar a $B$ como intermediario. En lenguaje matemático, a la flecha $h$ de $A$ a $C$ se le llama la composición de las flechas $f$ y $g$, la cual se escribe $g\circ f$. La tercera propiedad tiene que ver con la forma en la cual, bajo estas composiciones, se comportan las flechas que salen y entran de una misma persona. Si para una persona cualquiera $X$ escribimos $1_X$ como la flecha resultado de su capacidad para asombrarse a si misma, cuando tenemos una flecha $f$ de una persona $A$ a otra $B$, la composición $1_B\circ f$ o $f\circ 1_A$ es de nuevo $f$, pues en ambos casos  la composición indica la misma información que $f$, es decir, $A$ puede sorprender a $B$.
	
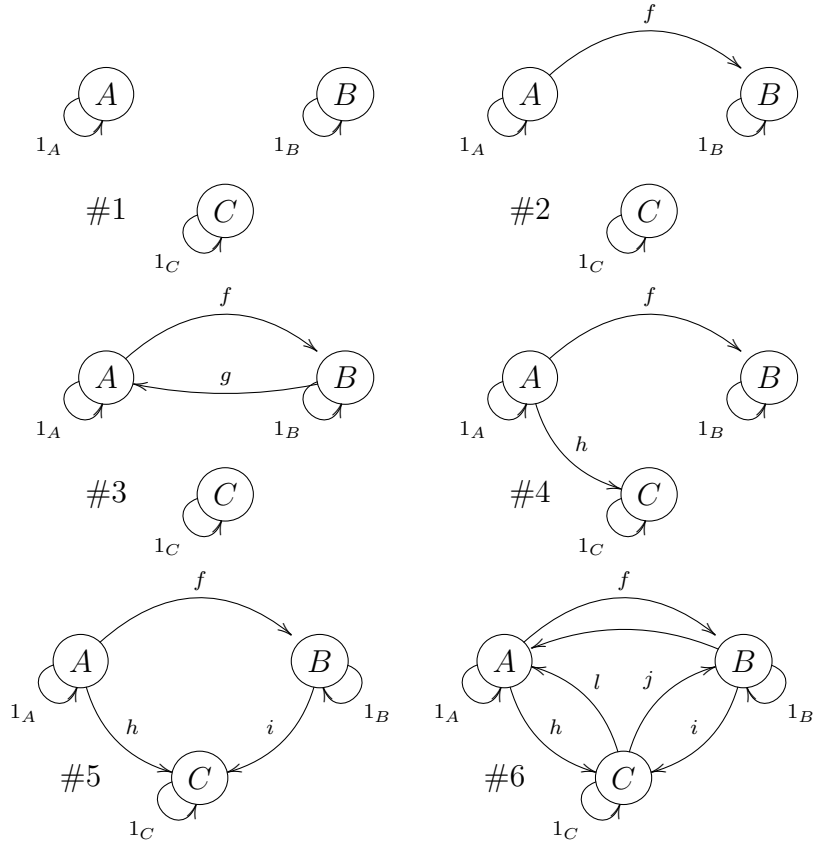
\begin{figure}[!]
\[\begin{matrix}
\entrymodifiers={++[o][F-]}
	\xymatrix@+0pc{
	A \ar@(l,d) []_{1_A} & *{}  &  B \ar@(l,d) []_{1_B}\\
		*{\# 1} & C \ar@(l,d) []_{1_C}&  *{}
	}
	&
	\entrymodifiers={++[o][F-]}
	\xymatrix@+0pc{
	A \ar@(l,d) []_{1_A} \ar@/^2pc/[rr]^f& *{}  &  B \ar@(l,d) []_{1_B}\\
		*{\# 2} & C \ar@(l,d) []_{1_C}&  *{}
	}\\ 
	\entrymodifiers={++[o][F-]}
	\xymatrix@+0pc{
	A \ar@(l,d) []_{1_A} \ar@/^2pc/[rr]^f& *{}  &  B \ar@(l,d) []_{1_B} \ar@/^0.5pc/[ll]_g\\
		*{\# 3} & C \ar@(l,d) []_{1_C}&  *{}
	}
	&
	\entrymodifiers={++[o][F-]}
	\xymatrix@+0pc{
	A \ar@(l,d) []_{1_A} \ar@/^2pc/[rr]^f \ar@/_1pc/[rd]^h& *{}  &  B \ar@(l,d) []_{1_B} \\
		*{\# 4} & C \ar@(l,d) []_{1_C}&  *{}
	}\\ 
	\entrymodifiers={++[o][F-]}
	\xymatrix@+0pc{
	A \ar@(l,d) []_{1_A} \ar@/^2pc/[rr]^f \ar@/_1pc/[rd]^h& *{}  &  B \ar@(r,d) []^{1_B} \ar@/^1pc/[dl]_i\\
		*{\# 5} & C \ar@(l,d) []_{1_C}&  *{}
	}
	&
	\entrymodifiers={++[o][F-]}
	\xymatrix@+0pc{
	A \ar@(l,d) []_{1_A} \ar@/^2pc/[rr]^f \ar@/_1pc/[rd]^h& *{}  &  B \ar@/_1pc/[ll]\ar@(r,d) []^{1_B} \ar@/^1pc/[dl]_i\\
		*{\# 6} & C \ar@(l,d) []_{1_C} \ar@/^1pc/[ru]^j \ar@/_1pc/[lu]_l  &  *{}
	}	
	\end{matrix}
	\]
\caption{Algunos tipos de categorías de tres objetos cuando solo puede haber una flecha de un objeto a otro. Tipo 1 representa la categoría con la mínima cantidad de flechas posibles, es decir, cada persona sólo es capáz de asombrarse a sí misma.
	Tipo 2, A asombra a B.
	Tipo 3, A asombra a B y a C.
	Tipo 4, A asombra a B y a C y B asombra a C.	
	Tipo 5, A asombra a B y B asombra a A.
	Tipo 6, tiene trazadas todas las posibles flechas entre estos objetos.}
\label{fig-cat-01}
\end{figure}

\FloatBarrier

	Notemos que el hecho de que $A$ pueda asombrar a $B$, no asegura que $B$ tenga la capacidad de generar la misma reacción en $A$, en otras palabras, la existencia de una flecha de $A$ a $B$ no implica la existencia de una flecha de $B$ a $A$. En la figura \ref{fig-cat-01} se representan diferentes versiones válidas de la categoría $\mathcal{P}$ con tres personas.


	Las situaciones que no son consideradas categorías son aquellas donde la composición de las flechas no se pueda realizar o no tenga sentido. Por ejemplo, si en una categoría como $\mathcal{P}$,	para cada par de personas $A$ y $B$, admitimos una flecha diferente de $A$ hasta $B$ por cada emoción que $A$ pueda provocar en $B$. Si $A$ puede asustar a $B$ y $B$ puede hacer sonreir a $C$, entonces ¿qué reacción puede provocar $A$ en $C$? Otro ejemplo de una situación donde $\mathcal{P}$ no sería considerada un categoría es cuando tomamos en cuenta el acto de contratar, pues no es cierto que si $A$ puede contratar a $B$ para un negocio y $B$ piensa que $C$ es acto para ser contratado, entonces sería falso que $A$ puede contratar a $C$, ya que los criterios de contratación entre $A$ y $B$ pueden variar. También se puede pensar en la situación donde se trazan las flechas de $\mathcal{P}$ a partir del acto de saludar, pues en este caso el resultado tampoco sería necesariamente una categoría, porque saludar implica un desconocimiento del otro, por lo tanto, saludarse a sí mismo no tedría sentido.

	Formalmente, una categoría $\mathcal{C}$ consiste de una colección de objetos $\text{Obj}(\mathcal{C})$ y una de flechas $\text{Fl}(\mathcal{C})$. Además, estas colecciones deberán cumplir con las siguientes condiciones:
	(1) Toda flecha $f$ de $\text{Fl}(\mathcal{C})$ inicia y termina en un objeto de $\text{Obj}(\mathcal{C})$.
	Cuando una flecha $f$ inicia en $A$ y termina en $B$, se escribe $f:A\to B$. 
	(2) Por cada objeto $X$ de $\text{Obj}(\mathcal{C})$, existe una flecha
	especial $1_X$ que inicia y termina en $X$ y se le llama flecha identidad de $X$.
	(3) Es posible componer flechas, es decir, por cada par de flechas $f:A\to B$ y $g:B\to C$,
	existe una regla que les asocia una única flecha $h:A\to C$, la cual se escribe $g\circ f:A\to C$.
	Además, el resultado de esta regla al componer varias flechas, no debe de estar 
	afectado por el orden de aplicación, es decir, si tenemos tres flechas
	$f:A\to B$, $g:B\to C$ y $h:C\to D$, aunque es posible componerlas de dos
	formas diferentes, iniciando con $f$ y $g$ para obtener $h\circ (g\circ f)$ o iniciando con
	$g$ y $h$ para obtener $(h\circ g)\circ f$, en ambos
	casos la flecha resultante de $A$ hasta $D$, debe ser la misma,
	así $h\circ (g\circ f)=(h\circ g)\circ f$.
	(4) Para cualquier flecha $f:A\to B$, las identidades $1_A$ y $1_B$ satisfacen
	$f = 1_B \circ f = f\circ 1_A$, es decir, no afectan el resultado de la composición.
	En una categoría $\mathcal{C}$, al conjunto de todas las flechas de $A$ hasta $B$,
	se le escribe $\text{Hom}(A,B)$, \cite{MacLane-1998}. 
	También, la palabra morfismo es usado como sinónimo de flecha.

	En matemáticas, las categorías surgen cuando se estudian estructuras con naturalezas similares.
	A continuación una lista con algunas de ellas.

	\textbf{La categoría de los conjuntos $Set$:}
	está formada por todos los conjuntos y las funciones entre ellos,
	es decir, $\text{Obj}(Set)$ son los conjuntos y $\text{Fl}(Set)$ son las funciones.

	\textbf{La categoría de los ordenes parciales $Pos$:}
   Los conjuntos que poseen una relaci\'on de orden parcial, digamos $(A,\leq)$ y $(B,\leq)$, 
	tambi\'en constituyen la colecci\'on de objetos de una categor\'ia denotada $Pos$, 
	donde las flechas $f\in \text{Hom}(A,B)$ son funciones mon\'otonas $f\:A\to B$.

 	\textbf{La categoría de los grupos $Grp$:}
	Un grupo $G$ es un conjunto junto con una operación binaria asociativa
	y un elemento especial $1_G$ llamado unidad del grupo, 
	el cual funciona como neutro de la operación, además de que
	por cada elemento del grupo existe otro con el cual, al
	aplicarle la operación binaria da como resultado $1_G$.
	Entre grupos hay funciones especiales que preservan la operación binarias,
	a una función de este tipo se le llama homomorfismo. 
	Los grupos junto con los homomorfismos forman una categoría $Grp$.

	\textbf{La categoría de los anillos $Rng$:}
	En un estilo similar al de un grupo, un anillo es un conjunto dotado de
	dos operaciones binarias $+$(suma) y $*$(producto). Por separado con cada
	operación, el conjunto debe se un grupo, pero además se deben de cumplir
	ser conmutativo bajo la primera operación y un algunas propiedades de
	distributividad entre las operaciones. Al igual que en grupos, las
	funciones importantes entre anillos son aquellas que preservan las operaciones y
	se llaman homomorfismos de anillos. Los anillos junto con sus homomorfismos
	forman una categoría denotada $Rng$.
		
	\textbf{La categoría de los espacios vectoriales $Vec$:} 
	Los espacios vectoriales junto con sus transformaciones lineales forman 
	uno de los primeros ejemplos de categorías presentados 
	por Saunders Maclane \cite{saunders-general-theory-of-natural-equivalences}.
	La categoría de espacios vectoriales se denota por $Vect$.

	\textbf{La categoría de los espacios topológicos $Top$:} 
	En el cálculo diferencial e integral de una o varias variables
	se estudian las principales propiedades de las funciones continuas.
	Dichas propiedades dependen también del dominio de las funciones.
	Así, el concepto de espacio topológico surge 
	a la hora de generalizar esta teoría a contextos más generales,
	Un espacio topológico es un conjunto donde se identifica una
	colección particular de sus subconjuntos, los cuales deben de
	satisfacer: (1) el conjunto principal y el conjunto vacío forman parte de esta colección.
	(2) El conjunto resultado de intersecar cualquier cantidad finita de estos
	subconjuntos, es de nuevo un subconjunto de la colección.
	(3) Cualquier conjunto resultado de una reunión de conjuntos de la colección
	es de nuevo un conjunto de la colección. 
	A la colección se le conoce como topología sobre el conjunto y a los conjuntos de la colección,
	se le llaman abiertos. 
	Una función continua entre dos espacios topológicos es toda aquella cuyo conjunto de
	preimágenes de cualquier abierto es un abierto. 
	La categoría de espacios topológico y funciones continuas se escribe $Top$.

\begin{figure}[h]
\[\entrymodifiers={++[o][F-]}
	\xymatrix@+0pc{
	A 
	\ar@(u,l) []_{1_A} 
	\ar@(l,dl) []_{A}  
	\ar@(dl,d) []_{AA}  
	\ar@(d,dr) []_{AAA} 
	  \ar@/^1pc/[rr]^{AM} 
	  \ar@/^3pc/[rrr]^{\,\,\,\,\,\,AMMAT}	 
	  \ar@/^5pc/[rrrr]^{AMMATTAS}
	& *{} 
	& 
	M \ar@(l,d) []_{1_M} 
	  \ar@/^1pc/[r]^{MAT\,\,\,\,\,\,\,\,\,\,\,\,\,\,\,\,} 
	& T \ar@(l,d) []_{1_T}  
	  \ar@/^1pc/[r]^{TAS\,\,\,\,\,\,\,\,\,\,\,\,} 
	& S \ar@(l,d) []_{1_T}  
	}
 \]
\caption{Algunas flechas de la categoría de palabras $Mots$.
	La flecha AA se puede obtener al componer la flecha A consigo misma.
	AAA es el resultado de componer A con AA. AM compuesto con MAT forma
	la flecha AMMAT, y esta misma compuesta con TAS, forma la flecha AMMATTAS.
	La flecha SSSS es ejemplo de una palabra que inicia y termina con S.
	$1_A$, $1_M$, $1_T$ y $1_S$, son las flechas identidad de cada objeto.}
\label{fig-cat-02}
\end{figure}
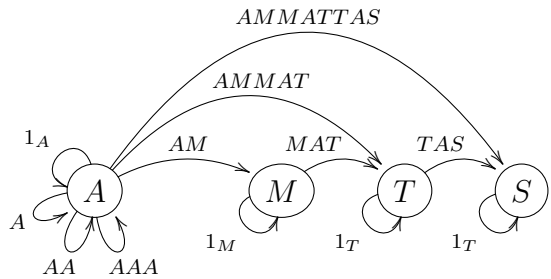

\textbf{La categoría de las palabras $Mots$:} En esta categoría los objetos
	son las letras del abecedario y cada palabra representa
	una flecha desde su letra inicial hasta su letra final.
	Las palabras de una sola letra, se ven como flechas que
	salen y llegan a la misma letra, aunque no se considerarían
	como morfismos identidad, la cual se puede ver como la palabra vacía. 
	La composición es el resultado
	de la concatenación de las palabras. 
	En $Mots$ las flechas podrían incluso representar palabras que sin significado
	en español, además de que tendríamos una infinidad de flechas diferentes entre objetos,
	incluso saliendo y llegando al mismo lugar, como lo son
	las flechas dadas por las palabras: \text{A, AA, AAA, AAAA},
	todas pertenecientes a $\text{Hom}(A,A)$.
	La figura \ref{fig-cat-02} es una representación de 
	algunas flechas en esta categoría.

\subsection{Morfismos especiales}
	Cuando dos estructuras de la misma naturaleza son indistingibles
	se dice que son isomorfas. Esto se refleja en teoría de categorías
	mediante la existencia entre los objetos de flechas especiales
	llamadas isomorfismos. La propiedad que las caracteriza
	es la siguiente: $f:A\to B$ es un isomorfismo si y solo si
	en la misma categoría existe $g:B\to A$ que satisface
	$f\circ g = 1_B$ y $g \circ f = 1_A$.

	En el caso de la categoría de conjuntos $Set$, como la propiedad
	que caracteriza a cada conjunto es sólo la cantidad de sus elementos,
	dos conjuntos son isomorfos cuando tienen la misma cantidad de elementos.
	Así, en $Set$
	un conjunto con cinco vasos es isomorfo a un conjunto de cinco platos.
	
\begin{figure}[h]
	\[
	\xymatrix{
	\text{\includegraphics[scale=0.5]{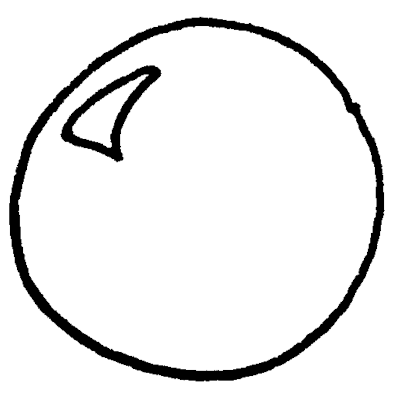}}\ar@<1ex>@/^2pc/@2{->}[r] & 
	\text{\includegraphics[scale=0.5]{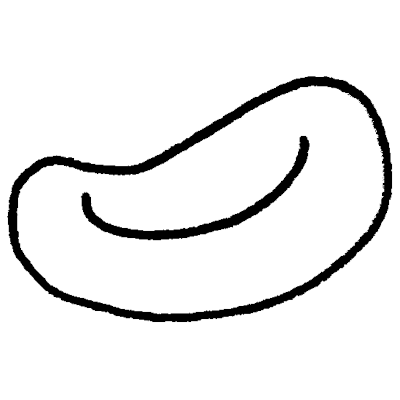}}\ar@<-2ex>@/^0.8pc/@2{->}[l] 
	}\]
\caption{Isomorfismo topológico en una pelota anti-estrés.}
\label{fig-cat-03}
\end{figure}

	Si cambiamos de categoría y nos vamos hacia la categoría de 
	espacios topológicos $Top$, los isomorfismos son un poco diferentes.		
	Tomemos un objeto cotidiano que podamos ver como un espacio topológico,
	por ejemplo una pelota anti-estrés. Cuando la sostenemos en la palma de la
	mano y la apretamos entre nuestros dedos, las deformaciones
	que sufre la pelota la transforman en un objeto con una forma diferente
	a la original pero no al grado de desnaturalizarla.
	Este tipo de acción es una representación de una aplicación continua
	entre dos espacios topológicos, la pelota en su forma original y la pelota bajo 
	la opresión de la mano, como se muestra en la figura \ref{fig-cat-03}. 
	Pero nuestra acción también representa
	un isomorfismo, ya que la pelota no ha sido desnaturalizada por la fuerza de nuestra mano.
	La naturaleza de este isomorfismo se manifiesta al dejar
	de ejercer fuerza sobre la pelota, pues automáticamente la pelota
	vuelve a su forma original usando otra aplicación continua,
	la cual sería inversa a nuestra fuerza original, aunque
	la manera en la que volvió a su estado original, no repite
	necesariamente nuestra acción al revés, es solo una forma de volver.

	Otros tipos de flechas importantes en una categoría son aquellas
	que generalizan la noción en la categoría $Set$ de subconjunto
	y proyección. Estas flechas son llamadas monomorfismos y epimorfismos,
	respectivamente. 
	
	En $Set$ un monomorfismo se da cuando un conjunto es un subconjunto de otro mayor,
	o al menos se puede identificar por un isormorfimo a un subconjunto de un conjunto.
	Lo interesante es que esta propiedad se puede codificar como una propiedad
	que habla sólo de flechas: $f$ es un monomorfismo si para cualquier
	par de flechas $g,h$, siempre que $f\circ g=f\circ h$, se satisface que $g=h$.
	En $Set$, los isomorfismos son exactamente las funciones inyectivas.
 
 	Por otro lado, el concepto de proyección entre conjuntos se da cuando 
        se tiene una función desde un conjunto hacia un subconjunto de él mismo,
        de tal manera que todos los elementos del subconjunto son usados,
	lo cual coincide con el concepto de función sobreyectiva.
	Al igual que el caso de los monomorfismos, esta propiedad se puede 
	expresar exclusivamente en términos de flechas y así, se puede exportar
	a cualquier categoría: se dice que $f$ es un epimorfismo si 
	para todo par de flechas $g$ y $h$ tales que $g\circ f=h\circ f$,
	se tiene que $g=h$.

 \subsection{Funtores entre categorías}
 
	Las categorías, al igual que los objetos que las forman, 
	se relacionan por entidades que se pueden
	comprender como flechas entre ellas, llamadas
	funtores. Estos se dividen en dos tipos, los covariantes
	y los contravariantes.
	
	Existe un funtor covariante $F$ de una categoría $\mathcal{C}$
	a otra categoría $\mathcal{D}$
	cuando se pueden crear dos tipos de aplicaciones, 
	una entre las colecciones de objetos 
	y otra entre las colecciones de flechas, de tal 
	manera que se respeten las propiedades de la composición y
	los morfismos unidad. Así, un funtor $F:\mathcal{C}\to \mathcal{D}$,
	asocia a cada objeto $X$ en $\mathcal{C}$, un objeto $F(X)$ en $\mathcal{D}$,
	y a cada flecha $f:A\to B$ de $\mathcal{C}$, una flecha $F(f):F(A)\to F(B)$ de $\mathcal{D}$.
	Además, se debe de cumplir que $F(1_X)=1_{F(X)}$ y $F(f\circ g)=F(f)\circ F(g)$.
	Cuando el funtor es contravariante, lo que cambia es el sentido
	de las flechas, es decir, si $f:A\to B$ en $\mathcal{C}$, 
	entonces $F(f)$ en $\mathcal{D}$ sería una flecha $F(f):F(B)\to F(A)$.
	También, en un funtor contravariante cambia el comportamiento de las composiciones,
	pues $F(f\circ g)=F(g)\circ F(f)$. 
	
	Si tomamos un categoría cualquiera $\mathcal{C}$ y cambiamos la
	dirección de cada una de sus flechas, entonces se obtiene
	lo que se llama la categoría opuesta de $\mathcal{C}$, la cual 
	se denota como $\mathcal{C}^{\text{op}}$. Este acto de invertir
	las flechas de $\mathcal{C}$, es un ejemplo de funtor contravariante 
	$F:\mathcal{C}\to \mathcal{C}^{\text{op}}$,
	el cual satisface $F(X)=X$ para todo objeto $X$ de $\mathcal{C}$,
	y $F(f):B\to A$ si $f:A\to B$. Una propiedad de las categorías opuestas
	es que los conceptos de monomorfismo y epimorfismo son duales,
	es decir, $f$ es un monomorfimos en $\mathcal{C}$ si y solo 
	si es un epimorfismo en $\mathcal{C}^{\text{op}}$.
	
	Un ejemplo trivial de funtor $F$, es aquel que va de una categoría
	en ella misma, $F:\mathcal{C}\to \mathcal{C}$, y asigna a cada
	objeto, un objeto $X$ fijo de $\mathcal{C}$, y a todo morfismo
	$f$, la identidad $1_X$.
 
	Otro tipo de funtores muy utilizado son los funtores de olvido. Con ellos
	se pueden describir 
	importantes construcciones matemáticas como los grupos libres o los espacios vectoriales generados.
	Un funtor de olvido toma objetos de una categoría y los desliga de algunas de sus propiedades,
	por ejemplo,
	un grupo $G$ es un conjunto con una operación binaria de ciertas condiciones,
	si privamos a $G$ de esta operación, es decir, si se olvida mencionar que $G$ tiene una operación
	binaria, $G$ es solo un conjunto.
	En cuando a los morfimos entre grupos,
	todos ellos son funciones entre los conjuntos subyacentes. 
	Esta asociación determina un funtor de olvido 
	$U:Grp \to Sets$. De forma similar, se pueden definir otros funtores de olvido
	$U:Rng\to Sets$, $U:Pos\to Sets$, $U:Vect\to Sets$ y $U:Top\to Sets$.

	Considere el funtor de olvido $U:Top \to Sets$, se quiere un funtor $F:Sets\to Top$ 
	que a cada conjunto le asigne un espacio topológico.
	Para definir $F$ se debe de pensar en construir un espacio topológico a partir de
	un conjunto cualquiera $X$, de tal manera que cualquier función entre conjuntos
	se convierta en una aplicación continua por la acción de $F$. 
	Esto se puede resolver, asignándole a cada conjunto $X$
	la topología formada por todos los posibles subconjuntos de $X$.
	Note que en este caso, si $X$ es un conjunto, el resultado de $U(F(X))$ es de nuevo el conjunto $X$,
	sin embargo, si $Y$ es un espacio topológico, $F(U(Y))$ no es necesariamente un espacio topológico
	igual a $Y$, pues la topología de $F(U(Y))$ es mayor a la de $Y$.

	Otro funtor topológico es 
	$F: Pos \to Top$, el cual transforma cualquier conjunto ordenado $X$ en un espacio topológico.
	Este funtor le asigna a $X$ la topología dada por la propiedad:
	$W$ es un abierto de $X$ si siempre que $x$ esté en $W$ y $y$ sea mayor que $x$, según el orden de $X$, entonces
	$y$ también está en $W$.
	
	Los funtores pueden aparecer en contextos menos matemáticos,
	como es el caso cuando analizamos las posturas de una persona para construir una categoría. 
		Esto se puede hacer de la siguiente manera:
		cada posible postura se considera un objeto de la categoría
		y existe una flecha de una postura $A$ a otra $B$
		por cada forma de mover el cuerpo para llevarlo de la
		postura $A$ a la postura $B$.
		Así, diferentes maneras de llegar a la misma postura
		representan diferentes flechas. También, podríamos decir que
		entre dos posturas hay una flecha si existe un video
		que inicie con una y termine en la otra.
			
 	Para dos personas diferentes, las categorías de movimientos asociadas
 	a cada una, en general, serían diferentes, ya que
 	la flexibilidad podría diferir entre ellas.
 		
\begin{figure}[h]
\includegraphics{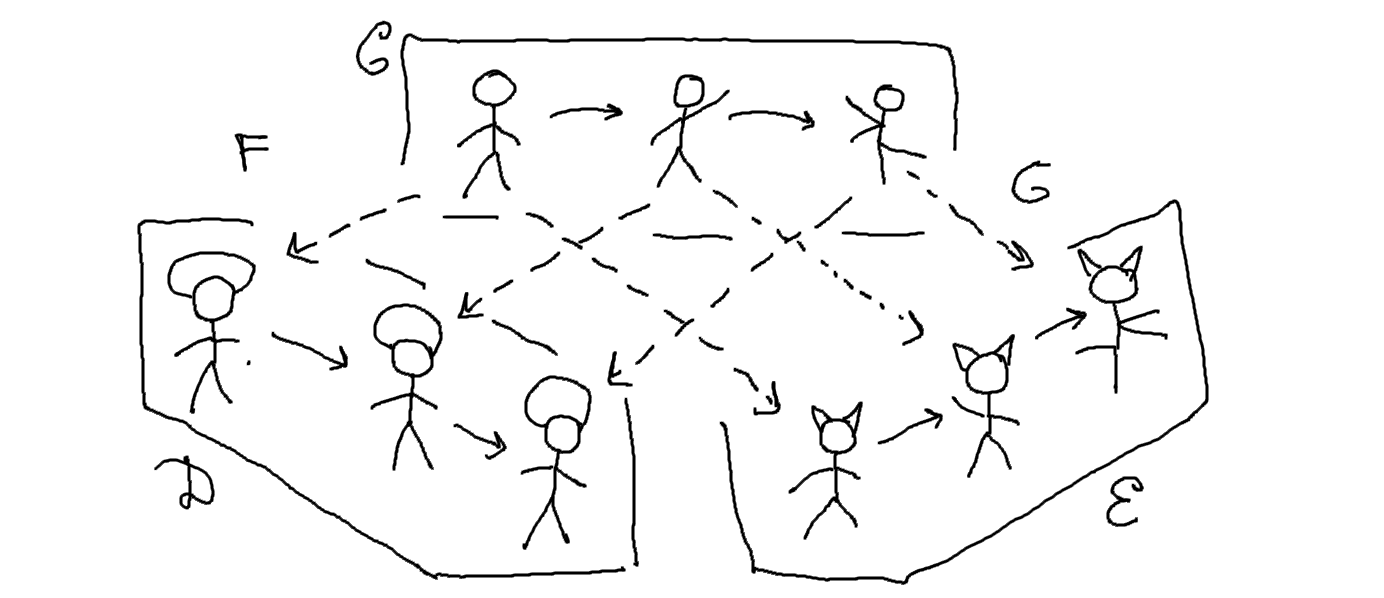}
\caption{Se representan tres distintas categorías de posturas C, D y E. De C a D está el
	funtor F de no hacer nada, donde cada movimiento se asocia a la misma postura.
	De C a E, las líneas punteadas indican el funtor de imitación en espejo.}
\label{fig-cat-04}
\end{figure}

 	Entre estas categorías surgen varios funtores que permiten
 	una clara representación gráfica. 	
		 Un ejemplo de estos funtores 
		esta dado por la imitación de los movimientos
		de una persona por otra persona.
		La acción del funtor imitación sería:
		a cada postura copia la misma postura y cada transición de
		postura se haría de la misma manera por parte de la otra.
		También podemos pensar en un funtor de imitación en espejo, 
		con el cual cada movimiento en una dirección es
		imitado en la dirección opuesta. 
		Un tercer ejemplo, es el funtor de no hacer nada,
		donde cada postura es correspondida por la misma postura estática,
		y cada cambio de postura es asignado a la acción de quedarse 
		estático en la misma postura. Este funtor y el de imitación en espejo
		están representados en la figura \ref{fig-cat-04}.

\subsection{Transformaciones naturales}

\begin{figure}[h]
\includegraphics{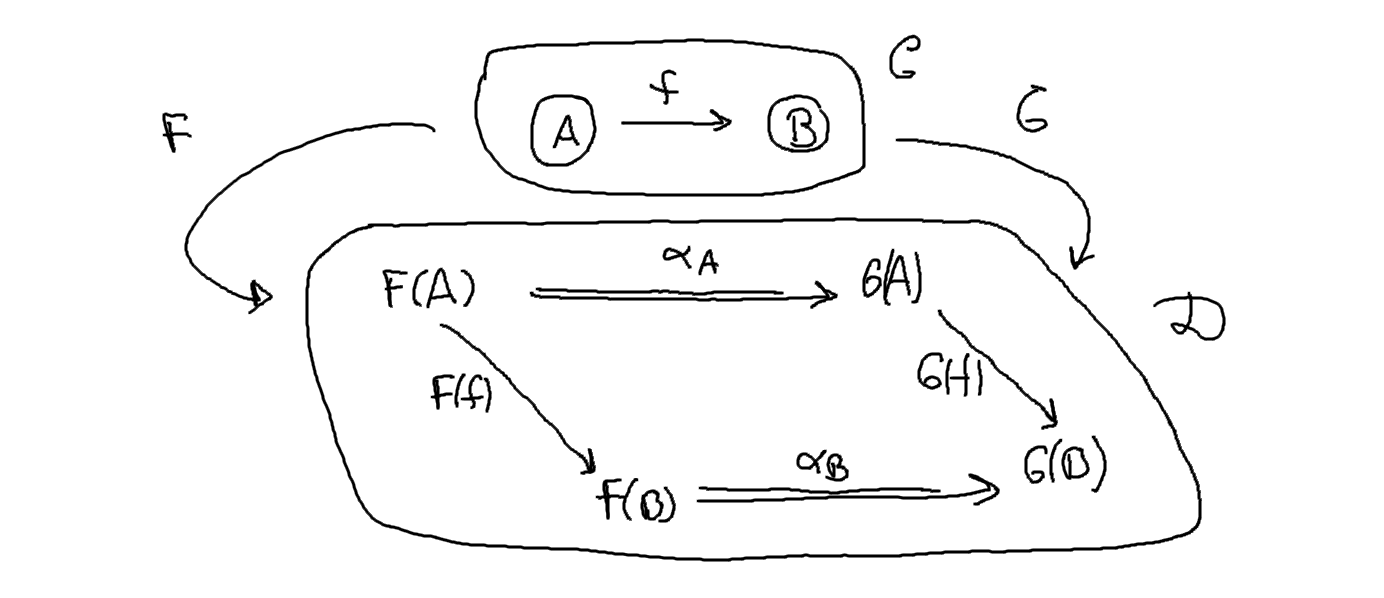}
\caption{F y G son dos funtores que van de la categoría C a la D, $\alpha$ representa
	una transformación natural de F a G, la cual genera el diagrama 
	mostrado cuando se le aplica a una flecha $f:A\to B$
	y que además satisface $\alpha_B \circ F(f)=G(f)\circ \alpha_A$. }
\label{fig-cat-05}
\end{figure}

	Las transformaciones naturales en teoría de categorías 
	son flechas entre funtores de tal manera que, los funtores
	entre dos categorías formen una categoría.
	En detalle, se definen de la siguiente manera:
	si $\mathcal{C}$ y $\mathcal{D}$ son dos categorías 
	y $F:\mathcal{C}\to \mathcal{D}$ y $G:\mathcal{C}\to \mathcal{D}$ 
	dos funtores entre ellas,
	una transformación natural $\alpha$ del funtor $F$ al funtor $G$,
	es una colección de flechas de la categoría $\mathcal{D}$
	de la forma $\alpha_X:F(X)\to G(X)$, una por cada objeto $X$ en $\mathcal{C}$.
	Además, estas flechas deben de satisfacer por cada flecha $f:A\to B$ de $\mathcal{C}$
	la condición de naturalidad $\alpha_B \circ F(f)=G(f)\circ \alpha_A$.	
	La figura \ref{fig-cat-05} ilustra esta igualdad en un diagrama.

\begin{figure}
\includegraphics{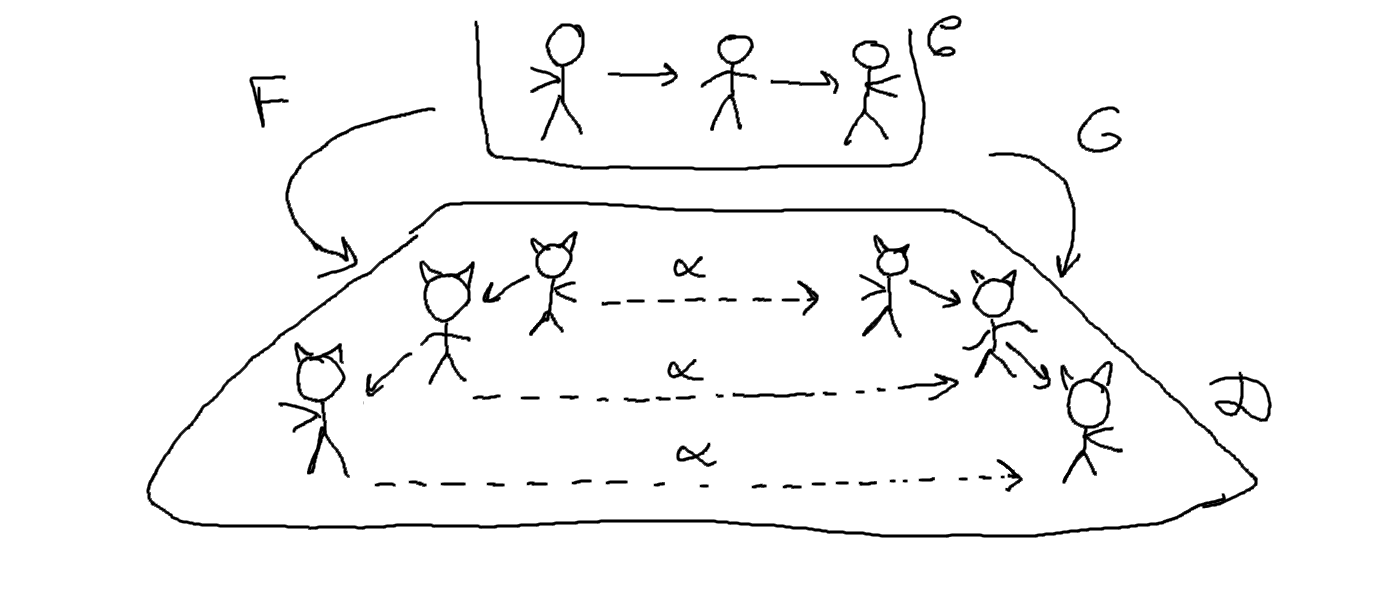}
\caption{Representación gráfica con las categorías de posturas de personas,
	de una transformación natural $\alpha$ entre dos funtores F y G, 
	como una coreografía entre dos copias de la misma persona 
	que guían sus movimientos con respecto a otra. }
\label{fig-cat-06}
\end{figure}

	Para representar gráficamente una transformación natural,
	retomemos el ejemplo de la categoría de las posturas de las personas. 
	Pensemos en $A$ y $B$ dos personas, y en sus categorías de posturas 
	$\mathcal{P}_A$ y $\mathcal{P}_B$. 
	Los funtores de $\mathcal{P}_A$ hacia $\mathcal{P}_B$, 
	representan los diferentes tipos de imitaciones o seguimientos que hace la
	persona $B$ de las posturas de $A$. 
	De esta manera, una transformación natural $\alpha$,
	entre dos funtores o seguimientos de $B$ de las posturas de $A$,	
	se puede ver como una coreografía de dos copias de la persona $B$
	al compás de los movimientos de $A$, 
	como se puede apreciar en la figura \ref{fig-cat-06}.

 \subsection{Lema de Yoneda}
 


	El lema de Yoneda \cite{MacLane-1998} enuncia una propiedad importante en el estudio
	de los funtores	que llegan a la categoría de los conjuntos.
	Su enunciado dice que para cualquier funtor  de la forma
	$F:\mathcal{C}\to Sets$ y cualquier objeto $X$ de la categoría $\mathcal{C}$,
	va a existir una relación uno a uno
	entre todas las posibles transformaciones naturales 
	del funtor $Hom_X$ al funtor $F$,
	y los elementos del conjunto $F(X)$. 
	Además, esta relación es natural, es decir, se
	conserva a través de cualquier morfismo $f:X\to Y$.
	En otras palabras, el lema Yoneda dice que los elementos de $F(X)$ codifican
	las transformaciones lineales entre $Hom_X$ y $F$.
	
	Para una representación del lema de Yoneda, pensemos en una
	persona $A$ y su categoría de posturas $\mathcal{C}$.
	Un funtor del tipo $F:\mathcal{C}\to Sets$ puede ser 
	uno que a cada postura le asigne un conjunto 
	de números que codifiquen algunas puntos importantes de las posturas,
	como las posiciones de lo brazos o piernas. Es claro que
	la complejidad del funtor dependerá que tan detallada queremos
	nuestra descripción de las posturas.
	Por ejemplo, a cada postura se le pueden
	asignar cuatro números de tres dígitos
	de tal manera que
	los códigos asociados a una postura se lean así: primer dígito 1 brazo, 0 pierna,
	segundo dígito 1 derecha, 0 izquierda, tercer dígito 1 arriba, 0 abajo.
	Si fijamos una postura $X$, para cualquier otra postura $Y$ de la persona, el conjunto
	$Hom_X(Y)=Hom(X,Y)$ son todos los movimientos posibles que puede
	hacer la persona para pasar de la postura $X$ a la $Y$.

	Ahora, una transformación natural $\alpha$ de $Hom_X$ a $F$
	es una correspondencia entre los movimientos de $X$ a $Y$,
	con los códigos de la postura $Y$, para cualquier $Y$.
	En particular, si $Y=X$ y se toma el movimiento constante de $X$ a $X$,
	denotado $1_X$, entonces $\alpha$ corresponde $1_X$ con un código de $X$,
	es decir, un elemento $x$ de $F(X)$.
	
\begin{figure}
\includegraphics{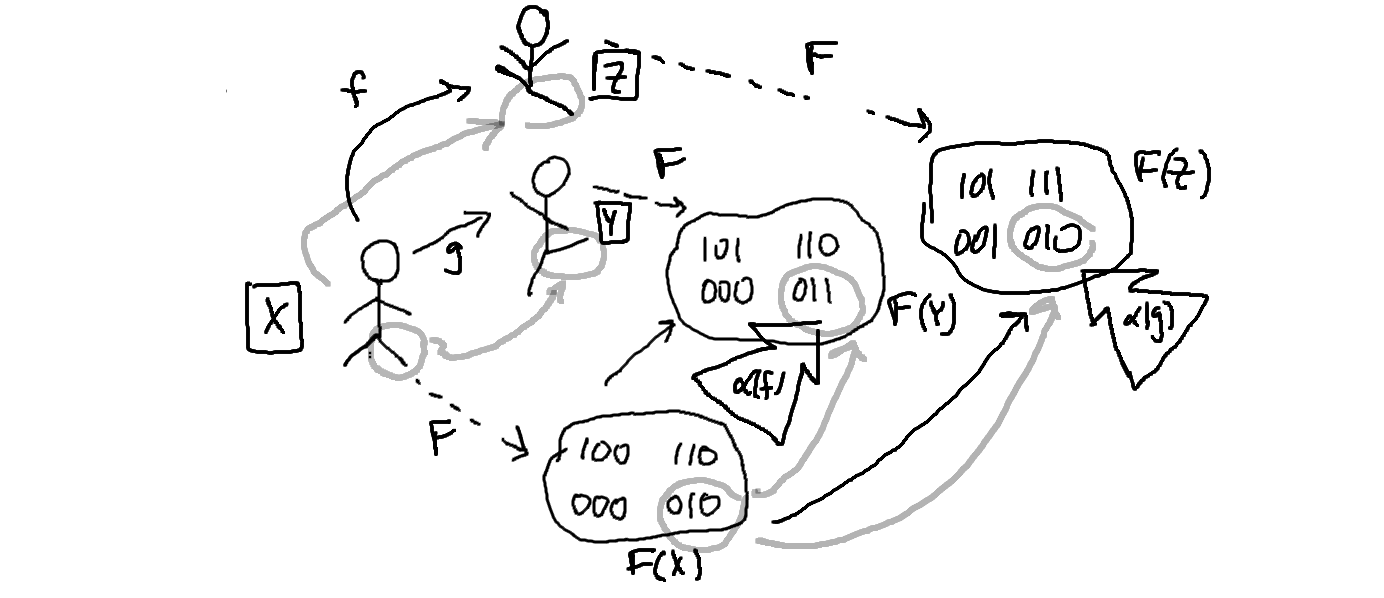}
\caption{Representación gráfica del lema de Yoneda. Se usan dos variaciones Y y Z
	de posturas de una persona a partir de una postura fija X. 
	El funtor F asigna a cada postura el conjunto de códigos. 
	La correspondencia expuesta en el lemma de Yoneda dice que
	cada código en F(X) de la
	postura base X tiene un único seguimiento a lo largo de los cambios de postura,
	dado por una única transformación natural $\alpha$ entre $Hom_X$ y F.}
\label{fig-cat-07}
\end{figure}

	Por otro lado, un código $x$ de la postura $X$, se refiera a cómo 
	se encuentra una parte del cuerpo de la persona. Cuando
	sucede un movimiento $f$ de $X$ a $Y$, si nos concentramos en la
	la parte descrita por el código $x$, este pasa a ser bajo $f$,
	el código $F(f)(x)$ en $F(Y)$, es decir, el código de la parte del cuerpo
	que representa $x$ pero con el valor de la postura $Y$, lo cual
	es una transformación natural de $Hom_X$ a $F$.
	En resumen, en el caso de las categorías de posturas, el lema de
	Yoneda nos dice que los números de una codificación $F$ de una postura $X$,
	representan todas las posibles imitaciones que puede hacer el 
	funtor codificador $F$ del funtor de movimientos $Hom_X$.
	La figura \ref{fig-cat-07} muestra a la correspondencia del lema de Yoneda en este contexto.
	
\section{Topos Elementales}
\label{sec-elemental-topoi}
Más allá del conocido universo de los conjuntos, donde la matemática clásica toma lugar,
se encuentran categorías cuyas estructuras llevan diferentes tipos de lógicas.
Algunas veces, estas estructuras son lo suficientemente robustas 
para ser consideradas universos matemáticos por sí mismos.
Es decir, en el contexto de la lógica matemática, se pueden construir proposiciones,
teoremas y demostraciones a partir de los objetos y flechas que pueblan la categoría.
A este tipo de categorías se les llama topos elementales \cite{johnstone2002sketches}.

La palabra topos, de origen griego y que significa lugar, es elegida por Alexander Grothendieck
al inicio de la década de los sesenta cuando usa este concepto como base de la teoría de cohomología etale,
desarrollada por él mismo \cite{bourbaki2006theorie}.



En su tesis, Lawvere ofrece una versión categórica de teorías algebraicas.
También sugiere que la categoría de categorías podría tomarse como
fundamento para las matemáticas y que los conjuntos podrían ser analizados en una
forma categórica. En los años siguientes, Lawvere trató de extender su análisis
y creó una primera versión categórica de teorías de primer orden
bajo el nombre de teorías elementales.
Luego, en 1969, en colaboración con Myles Tierney,
Lawvere introdujo la noción de topos elemental,
haciendo explícita la conexión entre la lógica de orden superior y la teorías de tipos \cite{Marquis2011-MARTHO-8}.

Los topos elementales fueron definidos con la idea de establecer 
condiciones mínimas para hacer posible la construcción de toda la maquinaria lógica matemática,
tomando como punto de partida el concepto más especializado de topos presentado por Grothendieck,
los cuales aparecen naturalmente en
el estudio de fibrados de espacios topológicos \cite{bourbaki2006theorie}.


\subsection{Funtores de subobjetos}

En aspectos más técnicos, 
el concepto más importante en teoría de categorías utilizado en la construcción de 
topos elementales es el funtor de subobjetos.
Este funtor contravariante y cuya categoría de llegada es $Sets$, 
se puede ver como una generalización 
del funtor asociado al conjunto de partes,
el cual asigna a cada conjunto, el conjunto de todos
sus subconjuntos.

Así, el funtor de subobjetos, a cada objeto de una categoría, asocia el conjunto
de todos los monomorfismos que llegan a él, eligiendo un solo representante
cuando los monomorfismos tienen objetos de salida isomorfos. 
Cuando esto se realiza sobre la categoría $Sets$, obtenemos el funtor de partes.

\begin{figure}
\includegraphics{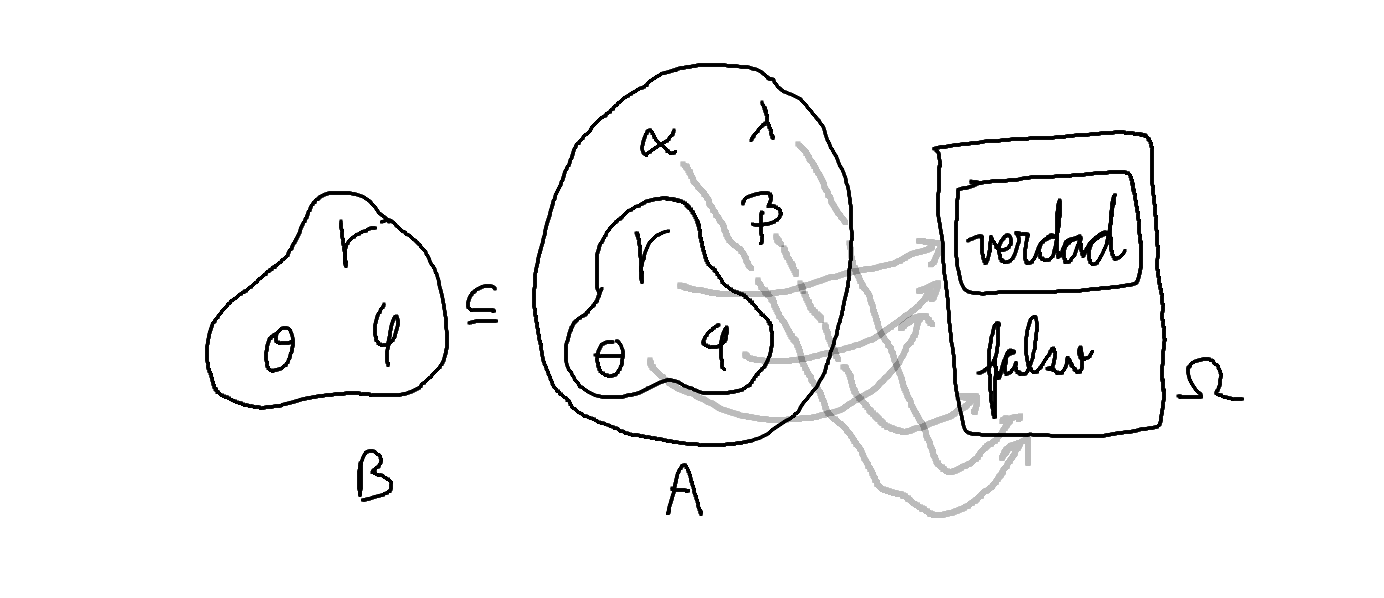}
\caption{B es subconjunto de A, el cual también se puede describir como
	el formado por los elementos de A asignados por una función al
	valor verdad.}
\label{fig-cat-08}
\end{figure}

El funtor de subobjetos posee una característica importante, es representable.
Para ilustrar esto, pasemos al contexto de $Sets$.
En $Sets$ el conjunto de partes de un conjunto tiene exactamente la
misma cantidad de elementos, que el conjunto de todas las funciones que
salen del conjunto y llegan al conjunto formado por los dos elementos ${\{\text{falso}, \text{verdadero} \}}$.
De esta forma, la información que necesitamos para
describir a los subconjuntos de un conjunto es el conjunto ${\{\text{falso}, \text{verdadero} \}}$
y un elemento particular de él, \emph{verdadero}. En efecto,
dada cualquier función desde nuestro conjunto hasta ${\{\text{falso}, \text{verdadero} \}}$,
basta ver quienes son asignados a \emph{verdadero}, para definir un subconjunto único,
como se muestra en la figura \ref{fig-cat-08}.
En un contexto más general, la representabilidad del funtor de subobjetos se refiere a la
capacidad de poderlo describir a partir de un conjunto fijo de subobjetos, llamado valores de verdad,  y un elemento
distinguible en él, llamado \emph{verdad}.

Prouté \cite{alain-logique-cat} 
define los topos elementales como categorías donde todos los límites finitos existen y
donde para cualquier objeto, el funtor que asigna a un objeto arbitrario el
conjunto de subobjetos del producto cartesiano de él con el objeto inicial,
es un funtor representable.

\subsection{Lógica en topos elementales}

En un topos elemental la existencia de límites implica posibilidad de trabajar con pullbacks y la existencia de un objeto terminal $1$, el cual tiene la característica de no modificar a un objeto cuando se hace el producto cartesiano con él. 
La segunda condición dice que el funtor de subobjetos es representable, pues es el 
mismo funtor que asigna a cada objeto el conjunto de subobjetos del producto cartesiano
de él con $1$. De este funtor en particular se obtiene el conjunto de valores de verdad
y el elemento verdad de su lógica interna.

Veamos un caso de un topos con más de dos valores de verdad. En la categoría de conjuntos,
un elemento solo puede tener dos estados con respecto a un conjunto dado,
no pertenece (falso) o pertenece al conjunto (verdad).
Cuando cambiamos a la categoría de grafos la situación cambia ligeramente.
Recuerde que un grafo dirigido está formado por nodos y aristas dirigidas o flechas, las cuales
van de un nodo a otro. Los diferentes estados de verdad para el topos formado por los grafos
y aplicaciones de grafos, se obtienen al considerar las diferentes situaciones que puede
tener una flecha con respecto a un grafo dado.

\begin{figure}
\includegraphics{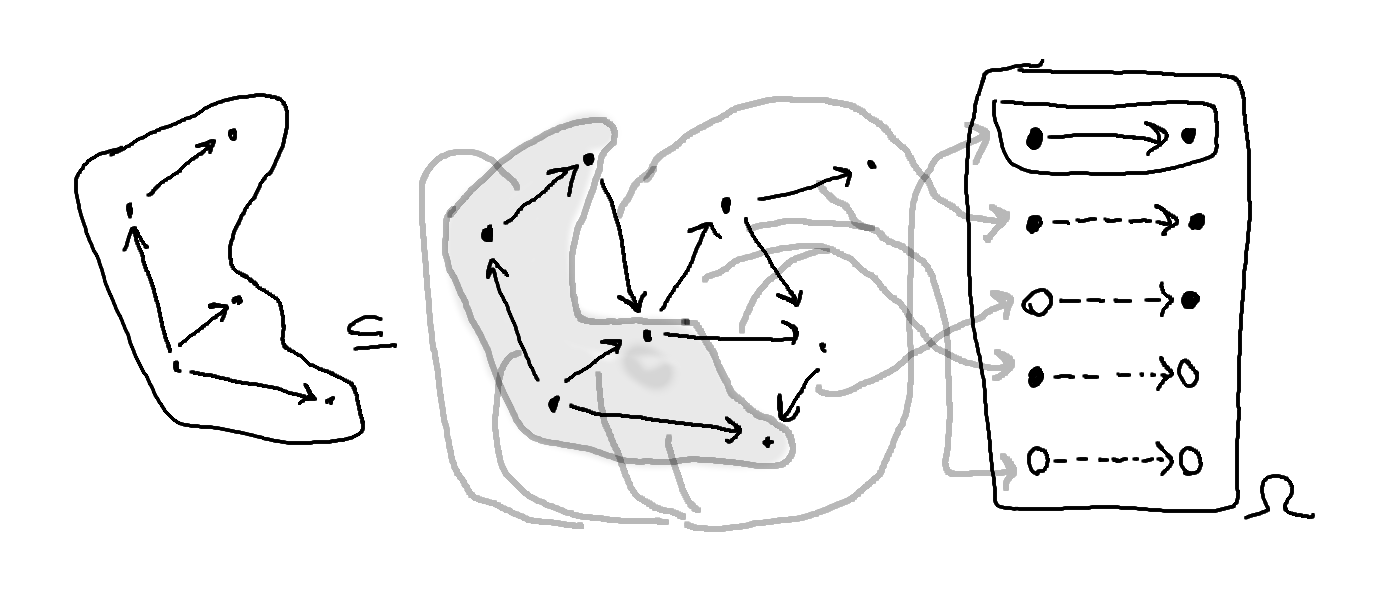}
\caption{Se representa un grafo y uno de sus subgrafos, además de una aplicación
	que hace corresponder cada flecha a su estado de verdad con respecto al subgrafo.}
\label{fig-cat-09}
\end{figure}

Sea $A$ un grafo y $B$ un subgrafo cualquiera de él. Para una flecha de $A$ tratemos de
responder a la pregunta ¿pertenece la flecha al subgrafo $B$? 
A diferencia de la categoría de los conjuntos, acá estamos en una situación donde 
aparecen diferentes estados entre lo que se puede considerar como verdadero y como falso.
De hecho, las diferentes situaciones son: 
(1) la flecha no pertenece completamente al subgrado (falso),
(2) solo su punto de partida pertenece al subgrafo,
(3) solo su punto de llegado pertenece al subgrafo,
(4) su punto de partida y llegada pertenecen al subgrafo, pero no su flecha,
y (5) toda la flecha pertenece al subgrado (verdad). En la figura
\ref{fig-cat-09} se hace una representación de los cinco estados
de una flecha con respecto a un subgrafo.

De esta forma, con el topos de grafos obtendríamos una lógica donde se tendrían
cinco estados de verdad, de los cuales tres son intermedios entre el falso y verdadero.
En la teoría de topos, se pueden construir topos con cantidades arbitrarias de 
estados de verdad entre falso y verdadero.

Prouté \cite{alain-logique-cat} da una clara analogía para comprender 
el desarrollo de la lógica interna en un topos elemental.
Las flechas que describen predicados y conectores lógicos pueden ser pensados
como el lenguaje de máquina de una computadora.
Para poder manipular todas estas flechas
se necesita un lenguaje de orden superior,
el cual es el lenguaje de Benabou-Mitchell o lenguaje interno.
Este proceso permite describir objetos en un topos
como lo hacemos en el universo de los conjuntos,
lo cual se puede entender como el aspecto sintáctico de la teoría.
Con respecto al aspecto semántico, se recurre a la semántica de Kripke-Joyal.
Esta es la manera en la cual la verdad de los predicados en la lógica interna
de un topos es interpretada como propiedades de sus flechas.

En una exploración más profunda de las propiedades de los topos elementales 
se podría plantear la pregunta sobre el ¿cuál es el grado de
parecido entre la lógica matemática ordinaria y la lógica en un topos?,
por ejemplo, uno se podría preguntar hasta qué punto el principio del tercer excluido 
se satisface en un topos, o también podríamos preguntarnos lo mismo con el axioma de elección \cite{alain-logique-cat}.
En el análisis lógico de un topos, cada uno de estas cuestiones es abordada desde dos tipos diferentes de enfoques,
una con respecto a un lenguaje interno y otra con respecto a un lenguaje externo de un topos,
lo cual se puede ver como un preámbulo a la complejidad y bastedad de este reciente campo matemático.
La persona intersada en continuar su lectura sobre este tema puede ver en \cite{landry2017categories}
un detallado catálogo de aplicaciones de la teoría de categorías.

\section*{Agradecimientos}
\addcontentsline{toc}{section}{Agradecimientos}
	
	Este trabajo ha sido financiado por la Universidad de Costa Rica a trav\'{e}s de la Vicerrectoría de Investigación. Espec\'{i}ficamente, el primer autor es financiado por la Escuela de Matem\'{a}tica {\sc em}at, y el segundo autor es financiado por la Escuela de Matem\'{a}tica {\sc em}at a trav\'{e}s del {\sc cimpa}, mediante el proyecto {\tt 821-C1-010}. 

        Los autores adem\'{a}s agradecen a la Escuela de Filosof\'{i}a {\sc ef-ucr}, por la invitaci\'{o}n a participar en el I Coloquio de L\'{o}gica, Epistemolog\'{i}a y Metodolog\'{i}a, evento en el que se gesta este trabajo. Especial agradecimiento al Dr. Lorenzo Boccafogli.

\bibliographystyle{amsplain}

\bibliography{principal-bibliography-philosophy}
\addcontentsline{toc}{section}{Referencias bibliogr\'{a}ficas}





\vspace{2ex}

\begin{center}
    $========================================$
\end{center}

\begin{flushright}
   {\em Jes\'{u}s E. S\'{a}nchez--Guevara}\\
  \small Escuela de Matem\'atica, Universidad de Costa Rica {\sc em}at--{\sc ucr}\\
  \small San Jos\'e 11501, Costa Rica\\
  \small e-mail: \texttt{jesus.sanchez\_g@ucr.ac.cr}\\
  \small {\sc orcid} \href{https://orcid.org/0000-0001-8993-1538}{\tt \color{black!29!blue} 0000-0001-8993-1538}

\vspace{2ex}

  {\em Ronald A. Z\'u\~niga--Rojas}\\
  \small Centro de Investigaci\'on en Matem\'atica Pura y Aplicada {\sc cimpa}\\
  \small Escuela de Matem\'atica, Universidad de Costa Rica {\sc em}at--{\sc ucr}\\
  \small San Jos\'e 11501, Costa Rica\\
  \small e-mail: \texttt{ronald.zunigarojas@ucr.ac.cr}\\
  \small {\sc orcid} \href{https://orcid.org/0000-0003-3402-2526}{\tt \color{black!29!blue} 0000-0003-3402-2526}
\end{flushright}

\begin{center}
    $========================================$
\end{center}

\end{document}